\documentclass[12pt]{amsart}
\usepackage{amssymb,latexsym,amsmath,amsthm}

\newtheorem{lemma}{Lemma}
\newtheorem{theorem}{Theorem}
\newtheorem{proposition}{Proposition}
\newtheorem{corollary}{Corollary}

\newcommand{\F}{\mathbb{F}}
\newcommand{\fq}{\mathbb{F}_q}
\newcommand{\rmv}[1]{}

{\theoremstyle{definition}
\newtheorem{example}{Example}
\newtheorem{remark}{Remark}
}

\def\BeginExample{\begin{example}}
\def\EndExample{\QED\end{example}}
\def\BeginRemark{\begin{remark}}
\def\EndRemark{\QED\end{remark}}

\def\boxitat#1#2{\vbox             
     {\hrule\hbox{\vrule\kern#1%
     \vbox{\kern #1\hbox{#2}\kern#1}\kern#1\vrule}\hrule}}
\def\enclose#1{\boxitat{0pt}{#1}}

\def\qedspace{\null}
\def\qedblack{\qedspace                    
  \lower.6pt\hbox{\vrule height7pt width 5pt}}
\def\qedwhite{\qedspace                    
  \lower.6pt\enclose{%
        \hbox{\vrule height7pt width 0pt\hskip6pt}}}
\def\QED{\qedwhite}

\def\beq{\begin{equation}}
\def\eeq{\end{equation}}

\usepackage{amssymb} 
\usepackage{graphicx}

\usepackage{epsfig}

\def\Fq{{\mathbb F}_q}

\def\CN{\mathcal{N}}
\def\alp{{i}}

\def\proof{\noindent {\bf Proof }}
\def\pr{{\mathbb P}}
\def\ex{{\mathbb E}}
\def\var{{\mathbb V}}

\def\F{\mathbb{F}}
\def\fq{\mathbb{F}_q}
\def\Fq{\mathbb{F}_q}

\begin{document}

\title[A probabilistic approach to value sets]{A probabilistic approach to value sets of polynomials over finite fields}

\author{Zhicheng Gao and Qiang Wang}
\thanks{
Research of authors was partially supported
by NSERC of Canada.}
\date{}
\address{School of Mathematics and Statistics\\
Carleton University\\
1125 Colonel By Drive\\
Ottawa, ON K1S 5B6\\
Canada}
\email{zgao@math.carleton.ca, wang@math.carleton.ca}



\begin{abstract}
In this paper we study the distribution of the size of the value set for a random polynomial with degree at most $q-1$ over a finite field $\fq$.
We obtain the exact probability distribution and show that the number of missing values tends to a normal distribution
as $q$ goes to infinity. We obtain these results through a study of a random $r$-th order cyclotomic mappings. A variation on the size of the union of some random sets is also considered.
\end{abstract}

\maketitle

\section{Introduction}

Let $\F_q$ be the finite field of $q$ elements with characteristic $p$.  Let $\gamma$ be a fixed primitive element of $\mathbb{F}_q$ throughout the paper.
 The {\it value set} of a polynomial $g$ over $\fq$ is the set $V_g$ of images  when we view $g$ as a mapping from $\fq$ to itself.
 Clearly $g$ is a {\it permutation polynomial (PP)} of $\fq$ if and only if the
cardinality $|V_g|$ of the value set $V_g$ is $q$.
Asymptotic formulas such as $|V_g| = \lambda(g) q + O(q^{1/2})$, where $\lambda(g)$ is a constant depending only on certain Galois groups associated to $g$,
can be found in Birch and Swinnerton-Dyer \cite{BW}  and Cohen \cite{Cohen:70}.
Later, Williams \cite{Williams} proved that almost all polynomials $g$ of degree $d$ satisfy  $\lambda(g) = 1 -\frac{1}{2!} + \frac{1}{3!} + \cdots + (-1)^{d-1} \frac{1}{d!}$.

There are also several results on explicit upper bound for $|V_g|$ if $g$ is not a PP over $\fq$; see for example \cite{GurWan,Tu,W1}.
 Perhaps the most well-known result is due to  Wan \cite{W1} who proved that if a polynomial $g$ of degree $d$ is not a PP then

\begin{equation}\label{Wan's bound}
 \vert V_g \vert \leq q -{q-1\over d}.
\end{equation}


On the other hand, it is easy to see that $|V_g|\ge \lceil q/d \rceil$ for any polynomial $g$ over $\fq$  with degree $d$. The polynomials achieving this lower bound are called   {\it minimal value set polynomials}. The classification of minimal value set polynomials over $\F_{p^k}$ with $k\leq 2$ can be found in  \cite{Carlitzetal:61,Mills:64}, and in \cite{BC} for all the minimal value set polynomials in $\fq[x]$ whose value set is a subfield of $\fq$. See  \cite{DasMul,WSC} for further results on lower bounds of $|V_g|$ and \cite{GomezMadden:88} for some classes of polynomials with small value sets.  More recently, algorithms and complexity in computing $|V_g|$ have been studied in \cite{CHW}.

We note that all of  these results mentioned above relate $|V_g|$ to the degree $d$ of $g$.  It is also well known that every polynomial $g$ over $\mathbb{F}_q$ such that $g(0) =b$
has the form $ax^rf(x^s)+b$ with some positive integers
 $r, s$ such that $s\mid q-1$. There are different ways to choose $r, s$ in the form
$ax^rf(x^s)+b$.   However,  in \cite{AGW:09},  the concept of the index of a polynomial was first introduced
and any  non-constant  polynomial $g\in\mathbb{F}_q[x]$ of
degree $\leq q-1$  can be
written {\it uniquely} as
$g(x) = a(x^rf(x^{(q-1)/\ell}))+b$
with  index $\ell$ defined below.
Namely, write
$$g(x)=a(x^n+a_{n-i_1} x^{n-i_1}+\cdots+a_{n-i_k} x^{n-i_k})+b,$$
where $a,~a_{n-i_j}\neq 0$, $j=1, \dots, k$.  The case that $k=0$ is trivial. Thus, we shall assume that
$k\geq 1$. Write $n-i_k=r$, the vanishing order of $x$ at $0$ (i.e., the lowest degree of $x$ in $g(x)-b$ is $r$).
Then $g(x)=a\left(x^r f(x^{(q-1)/\ell}) \right)+b,$ where $f(x)=
x^{e_0}+a_{n-i_1} x^{e_1}+\cdots+ a_{n-i_{k-1}}x^{e_{k-1}} + a_{r}
$,
 $$\ell=\frac{q-1}{\gcd(n-r,n-r-i_1,\dots, n-r-i_{k-1}, q-1)} := \frac{q-1}{s},$$
and $\gcd(e_0, e_1, \dots, e_{k-1}, \ell)=1 .$ The integer $\ell=\frac{q-1}{s}$ is called the
  {\it index} of $g(x)$.  From the above definition of index $\ell$, one can see that the greatest common divisor condition makes $\ell$ minimal among those possible choices.

Clearly, the study of the value set of $g$ over $\fq$ is equivalent to studying the value set $x^rf(x^{(q-1)/\ell})$ over $\fq$ with index $\ell$.
Recently Mullen, Wan and Wang \cite{MWW2} used an index approach to study the upper bound of the value set for any polynomial which is not a PP. They proved that if $g$ is not a PP then
\begin{equation}
|V_g| \leq q- \frac{q-1}{\ell}.
\end{equation}
This result improves Wan's result when the index $\ell$ of a polynomial is strictly smaller
than the degree $d$.  We note that the index $\ell$ of a polynomial is always smaller than the degree $d$ as long as $\ell \leq \sqrt{q} -1$.

The above result is obtained through a study of cyclotomic mapping polynomials which were studied earlier in  \cite{Evans:92,NW:05,Wang}. The index of a polynomial  is closely related to the concept of the least index of a cyclotomic mapping polynomial.  Recall that $\gamma$ is a fixed primitive element of $\mathbb{F}_q$. Let $\ell \mid q-1$ and the set of all nonzero $\ell$-th powers be $C_0$. Then $C_0$ is a subgroup of $\mathbb{F}_q^*$ of
index $\ell$. The elements of the factor group $\mathbb{F}_q^*/C_0$
are the {\it cyclotomic cosets}
$$ C_i := \gamma^i C_0, \  \  \  i = 0, 1, \cdots, \ell-1.$$

For any $a_0, a_1, \cdots, a_{\ell-1} \in \mathbb{F}_q$ and a positive integer $r$,   the  {\it $r$-th order cyclotomic mapping $f^{r}_{a_0, a_1, \cdots,
a_{\ell-1}}$ of index $\ell$ }  from $\mathbb{F}_q$ to itself  (see Niederreiter and Winterhof in \cite{NW:05} for $r=1$ or  Wang \cite{Wang}) is defined by
\begin{equation}\label{CycloMappingDef}
f^{r}_{a_0, a_1, \cdots,
a_{\ell-1}} (x) =
\left\{
\begin{array}{ll}
0, &   \mbox{if} ~ x=0; \\
a_i x^{r},   &   \mbox{if} ~x \in C_i, ~ 0\leq i \leq \ell-1. \\
\end{array}
\right.
\end{equation}

It is shown that  $r$-th order cyclotomic mappings of index $\ell$ produce the polynomials of the form  $x^r f(x^s)$ where $s =\frac{q-1}{\ell}$. Indeed,  the polynomial presentation is given by
\begin{equation*}
g(x) =\frac{1}{\ell} \sum_{i=0}^{\ell-1} a_i x^{r} \sum_{j=0}^{\ell-1} \zeta^{-ji}x^{js}, 
\end{equation*}
where $\zeta =\gamma^s$ is a fixed  primitive $\ell$-th root. On the other hand, as we mentioned earlier, each polynomial $f(x)$ such that $f(0)=0$ with index $\ell$ can be written as $x^r f(x^{(q-1)/\ell})$, which is an $r$-th order cyclotomic mapping with the least index $\ell$ such that $a_i = f(\zeta^i)$ for $i=0, \ldots, \ell-1$.

In this paper, we are interested in the probability distribution of the value set size of a random $r$-th order cyclotomic mapping polynomial for any given index $\ell$ and any positive integer $r$, as defined in Equation~(\ref{CycloMappingDef}). Thus this enables us to derive the probability distribution of the size of value set of a random polynomial of degree $d \leq q-1$ over a finite field.  In Section~\ref{Methodology} we outline our method
and a crucial result on normal distribution which is used in this paper. Essentially, we are interested in the distribution of the  size of the union of subsets, namely,  the distribution of the random variable $X_{t \ell}=|\cup_{j=1}^\ell A_j|$ where $A_j = g(C_{j-1})$ for $j=1, \ldots, \ell$ and $t=(r, s)$.  In Section~\ref{random mappings}, we first consider a simplified model such that none of $a_i$'s in Equation~(\ref{CycloMappingDef}) is zero.  Hence $a_i$'s are chosen independently at random from $\fq^*$. This means that the zero is not contained in any one of the subsets $A_j$'s. In particular, in Theorem~\ref{Thm1} we obtain the distribution of the number of missing values for a random $r$-th order cyclotomic mapping $f^{r}_{a_0, a_1, \cdots,
a_{\ell-1}} $ such that none of $a_i$'s is zero. Moreover, in Theorem~\ref{normaldistribution},  we show that this distribution is asymptotically normal.

In Section~\ref{random polynomials}, we study any random $r$-th order cyclotomic mapping polynomial by  choosing $a_i$'s in Equation~(\ref{CycloMappingDef}) independently at random from $\fq$. The probability distribution of value set size is given in Theorem~\ref{valuesetexclude0}. 
As a consequence, for $\ell = q-1$, we obtain the exact probability distribution of the value set size of a random polynomial over $\fq$ with degree at most $q-1$.
 In particular, we have the following corollaries to Theorem~\ref{valuesetexclude0}.

\begin{corollary}\label{Smallwithzero} Let $g(x)$ be a random polynomial of degree at most $q-1$ over $\fq$ with $g(0) =0$. Then
$$\pr(|V_g|= k+1)= \ {q-1\choose k}  \sum_{j=0}^{k} (-1)^{k-j} {k\choose j}\left(\frac{1+j}{q}\right)^{q-1}. $$
Consequently, for $k=o(q)$, we have
$$\pr(|V_g|= k+1)\sim \frac{1}{k!}(q-1)^k\left(\frac{k+1}{q}\right)^{q-1}.
$$
\end{corollary}

This proves that if $k > 1$ is small compared to $q$ then the number of polynomials over $\fq$ with degree less than or equal to $q-1$ such that the value set size is $k$ is always exponential in $q$. Moreover, we have

\begin{corollary}\label{randomnormaldistribution}
Let $g(x)$ be any random polynomial of degree at most $q-1$ over finite field $\fq$ with $g(0) =0$. Let $Y_q=q-|g(\fq)|$ denote the number of missing nonzero values in the value set of $g$. Let $\mu_{q} = q/e$ and $\sigma_q^2 = (e^{-1} - 2e^{-2})q$. Then the distribution of $(Y_q-\mu_q)/\sigma_q$ tends to the standard normal, as $q\to \infty$.
\end{corollary}

Finally in Section~\ref{union of random sets} we study a variation of our model used in Section~\ref{random mappings}. We consider the case when each subset $A_i$ is  chosen uniformly at random from all $m_i$-subsets of
a given $n$-set for $i=1, \ldots, \ell$. This extends a result by
Barot and Pe\~{n}a \cite{BP} and David \cite{Da}. We  also show that the size of the complement of $\cup_{i=1}^{\ell}A_i$ is asymptotically normal.

\section{Methodology}\label{Methodology}

In \cite{MWW}, we obtained the following formula for the cardinality of the value set for an arbitrary polynomial.

\begin{proposition}[Proposition 2.3 in \cite{MWW}]\label{valuesetCyclo}
Let $g(x) = a x^rf(x^s) + b $ ($a \neq 0$) be any polynomial over
$\fq$ with index $\ell = \frac{q-1}{s}$ and let $\gcd (r, s)=t$.  Let $\gamma$ be a fixed primitive element of $\fq$.
Then
\[|V_g| = c \frac{s}{t}  + 1,~or~ |V_g| = (c-1) \frac{s}{t}+1, \]
where  $c = |\{ ( \gamma^{ir} f(\gamma^{s i}) )^{s/t} \mid i = 0, \ldots, \ell-1 \}|$.
\end{proposition}

As discussed earlier, it is sufficient to assume that  $a=1$ and $b=0$ in Proposition~\ref{valuesetCyclo}. That is, we can view $g(x)$ as a $r$-th order cyclotomic mapping polynomial with the least index $\ell$.   In this case,
we have $g(x) = a_i x^r$ when $x\in C_i$, where $a_i = f(\gamma^{s i})$ for $i=0, \ldots, \ell-1$.
Recall that $C_0$ is the subgroup of $\fq^*$ consisting of all the $\ell$-th powers of $\fq^*$ and we let $T_0$ be the subgroup of $\fq^*$ consisting of all the $t\ell$-th powers. Hence
$T_i$ with $0\leq i \leq t\ell-1$ give all the cyclotomic cosets of index $t\ell$.
We also note that $x^r$ maps $C_0$ onto $T_0$ which contains  $\frac{s}{t}$ distinct elements.  So $x^r$ maps each coset $C_i = \gamma^i C_0$ onto $\gamma^{ir} T_0$.  Therefore $g$ maps $C_i$ onto  $ \gamma^{ir} f(\gamma^{s i})  T_0$, which could be either the set  $\{0\}$ (if $a_i = f(\gamma^{s i}) =0$) or one of the nonzero cyclotomic cosets of index $t\ell$. We observe that $c$ is the number of  distinct cyclotomic cosets of the form $\gamma^{ir}  f(\gamma^{s i}) T_0$, possibly along with the subset $\{0\}$ if one of $a_i$'s is zero. Hence we have  $|V_g| = c \frac{s}{t}  + 1$ or $(c-1) \frac{s}{t}+1$, the latter happens when  some of $a_i$'s in $g(x) = a_i x^r$ equal $0$.

Therefore the value set problem  for a random $r$-th order cyclotomic mapping polynomial (or random polynomial) $g$, essentially requires us to study the number $c$ in Proposition~\ref{valuesetCyclo}, the size of union of some cyclotomic cosets and possibly the subset $\{0\}$ if $a_i$'s take zero. More specifically, for $0\leq i \leq \ell-1$, each $C_i$ is mapped to $A_{i+1} =g(C_i)$ which is one of $T_0, \ldots, T_{t\ell -1}$ or $\{0\}$. Then $c$ is the number of distinct $A_j$'s  ($1\leq j\leq \ell$) and the value set size is either $c\frac{s}{t} + 1$ or $(c-1)\frac{s}{t} +1$.
More generally, we are interested in the distribution of the random variable $X_{t \ell}=|\cup_{j=1}^\ell A_j|$, while $A_j$ are  chosen independently according to a given distribution depending on that $a_i$ is chosen independently at random from $\fq$. Similar problems have been studied
in \cite{BP,Da} when each $A_j$ is chosen uniformly at random from all $k$-subsets of a given $n$-set.

Let $n=t\ell$ and let $D_0=\{0\}$ and $D_{j} = T_{j-1}$ for $1\leq j\leq t\ell -1$.  Let $Y_{n}$ be the the number of $D_1, \ldots, D_{n-1}$ which are not in $\cup_{j=1}^\ell A_j$ for a random $r$-th order cyclotomic mapping polynomial with index $\ell$ such that $(r, s)=t$.   We will derive exact probability distributions of $Y_n$ and show that they are asymptotically normal. Throughout the paper, we shall use $(Y_n)_k$ to denote the falling factorial $Y_n(Y_n-1)(Y_n-2)\cdots(Y_n-k+1)$. We use $\pr$, $\ex$, $\var$ to denote the probability, expectation, and variance of a random variable, respectively.

The main tool for deriving the probability distribution of $Y_n$ in the paper is through the sieve method and  the falling factorial moments.  Let
$B_1, \ldots, B_{n}$ be $n$ events in a probability space and $\CN =\{1, \ldots, n\}$.
We note that $\pr(Y_n=k)$ is the probability that exactly $k$ of the $B_j$ occur.  Define
$$S_h=\sum_{J\subset \CN,|J|=h}\pr(\cap_{j\in J}B_j).$$
Then the well-known sieve formula (See e.g. \cite[Theorm~10]{Bo}) gives
\begin{eqnarray}
\pr(Y_n=k)&=&\sum_{h=k}^n(-1)^{h-k}{h\choose k}S_h,\label{pk}\\
S_k&=&\sum_{h=k}^n{h\choose k}\pr(Y_n=h).\label{Sk}
\end{eqnarray}
Hence
\begin{equation}\label{factorial}
\ex((Y_n)_k)=k!S_k.
\end{equation}

We also  need the following
result \cite[Theorem~1]{GW} in order to show that $Y_n$ is asymptotically normal.

\begin{lemma}\label{normal}
Let $s_n>-\mu_n^{-1}$ and
$$\sigma_n=\sqrt{\mu_n+\mu_n^2s_n},$$
where $\mu_n\to \infty$ as $n\to \infty$.
Suppose that $$\mu_n=o(\sigma_n^3),$$ and a sequence $Y_n$ of nonnegative random variables satisfies
$$\ex((Y_n)_k)\sim \mu_n^k\exp\left(\frac{k^2s_n}{2}\right),$$
uniformly for all integers $k$ in the range $c\mu_n/\sigma_n\le k\le c'\mu_n/\sigma_n$ for some constants
$c'>c>0$. Then $(Y_n-\mu_n)/\sigma_n$ tends in distribution to the standard normal as $n\to \infty$.
\end{lemma}

\section{Sizes of value sets of cyclotomic mapping polynomials with nonzero branches} \label{random mappings}

In this section, we study of the value set size of a random $r$-th order cyclotomic mapping polynomial with nonzero branches (i.e., none of $a_i$'s is zero).  This means that we   choose $a_i$ in (\ref{CycloMappingDef}) independently at random from $\fq^*$ and it leads to the following model.
Let $n=t\ell$ and $D_1,D_2,\ldots D_{t\ell}$ be pairwise disjoint subsets of $\fq^*$
such that $|D_i|=s/t$ for all $1\le i\le t\ell$.  Because $a_i$ is chosen independently at random from $\fq^*$,  this means that   $A_1,A_2,\ldots, A_{\ell}$ are chosen independently and uniformly at random from $\{D_1,D_2,\ldots, D_{t\ell}\}$. We are interested in the distribution of $X_n=|\cup_{j=1}^{\ell}A_j|$. This is closely related to the distributions of $\ell$ labeled balls into $t\ell$ labeled boxes. Let $Y_n$ be the number of empty boxes in a random distribution of $\ell$ labeled balls
into $t\ell$ labeled boxes. We note $X_n=n-(s/t)Y_n$. We first prove the following
\begin{theorem}\label{Thm1}
Let $Y_n$ be the number of empty boxes in a random distribution of $\ell$ labeled balls into $n = t\ell$ labeled boxes.
We have
\begin{eqnarray*}
\ex((Y_n)_k)&=&(t\ell)_k\left(\frac{t\ell-k}{t\ell}\right)^{\ell},\\
\pr(Y_n=k)&=&{t\ell\choose k}\frac{1}{(t\ell)^{\ell}}\sum_{j=1}^{t\ell-k}(-1)^{t\ell-k-j}{t\ell-k\choose j}j^{\ell}.
\end{eqnarray*}
\end{theorem}
\proof
Let $\CN = \{1, 2, \ldots, t\ell\}$.  Let $B_j$ be the event that box $j$ is empty.
We have
$$S_k= \sum_{J\subset \CN,|J|=k}\pr(\cap_{j\in J}B_j) = {t\ell\choose k}\pr(\cap_{j=1}^k B_j)={t\ell\choose k}\left(\frac{t\ell -k}{t\ell}\right)^{\ell}.$$
It follows from (\ref{factorial}) that
$$\ex((Y_n)_k)=(t\ell)_k\left(\frac{t\ell-k}{t\ell}\right)^{\ell}.$$

Using Equation~(\ref{pk}), we obtain
\begin{eqnarray*}
\pr(Y_n=k)&=&\sum_{h=k}^{n} (-1)^{h-k} { h \choose k} S_h \\
&=& \sum_{h=k}^{t\ell} (-1)^{h-k} { h \choose k}  {t\ell\choose h}\left(\frac{t\ell -h}{t\ell}\right)^{\ell} \\
&=&  {t\ell\choose k} \sum_{h=k}^{t\ell} (-1)^{h-k} {t\ell-k\choose h-k} \left(\frac{t\ell -h}{t\ell}\right)^{\ell} \\
&=& {t\ell\choose k} \frac{1}{(t\ell)^\ell} \sum_{j=0}^{t\ell-k} (-1)^{t\ell -j-k} {t\ell-k\choose j} j^{\ell}\\
\end{eqnarray*}
\qed

Next we obtain

\begin{theorem}\label{normaldistribution}
Suppose $t=o(\ell^{1/5})$ as $n = t \ell \to \infty$. Define $$\mu_n=te^{-1/t}\ell,~\sigma_n^2=te^{-2/t}(e^{1/t}-1-1/t)\ell.$$
Then the distribution of $(Y_n-\mu_n)/\sigma_n$ tends to the standard normal, as $n\to \infty$.
\end{theorem}
\proof From Theorem~\ref{Thm1}, we have, as $\ell\to \infty$,
\begin{eqnarray}
\ex(Y_n)&=&t\ell\left(1-\frac{1}{t\ell}\right)^{\ell}\nonumber\\
&=&t\ell\exp\left(\ell\ln\left(1-\frac{1}{t\ell}\right)\right)\nonumber\\
&=&t\ell\exp\left(-\ell\left(\frac{1}{t\ell}+\frac{1}{2t^2\ell^2}+O\left(\frac{1}{(t\ell)^{3}}\right)\right)\right)\nonumber\\
&=&t\ell\exp\left(-\frac{1}{t}-\frac{1}{2t^2\ell}+O\left(\frac{1}{t^3\ell^2}\right)\right)\nonumber\\
&\sim&\mu_n, \label{mean} \nonumber \\
\end{eqnarray}
\begin{eqnarray}
\var(Y_n)&=&\ex(Y_n(Y_n-1))+\ex(Y_n)-(\ex(Y_n))^2\nonumber\\
&=&(t\ell)(t\ell-1)\left(1-\frac{2}{t\ell}\right)^{\ell} + t\ell \left(1-\frac{1}{t\ell}\right)^{\ell} - (t\ell)^2 \left(1-\frac{1}{t\ell}\right)^{2\ell} \nonumber\\
&=&(t\ell)(t\ell-1)\exp\left(\ell\ln\left(1-\frac{2}{t\ell}\right)\right)+ t\ell\exp\left(\ell\ln\left(1-\frac{1}{t\ell}\right)\right) \nonumber\\
&&-(t\ell)^2\exp\left(2\ell\ln\left(1-\frac{1}{t\ell}\right)\right) \nonumber\\
&=& (t\ell)(t\ell-1)\exp\left(-\ell\left(\frac{2}{t\ell} + \frac{2}{(t\ell)^2}+ O\left( \left(\frac{2}{t\ell}\right)^3 \right)\right)\right) \nonumber \\
&& + t\ell\exp\left(-\ell\left(\frac{1}{t\ell} + \frac{1}{2(t\ell)^2}+ O\left( \left(\frac{1}{t\ell}\right)^3 \right)\right)\right) \nonumber \\
&&  -(t\ell)^2\exp\left(-2\ell\left(\frac{1}{t\ell} + \frac{1}{2(t\ell)^2}+ O\left( \left(\frac{1}{t\ell}\right)^3 \right)\right)\right) \nonumber\\
&=& (t\ell)(t\ell-1)\exp\left(-\ell\left(\frac{2}{t\ell}\right)\right)\exp\left(-\frac{2\ell}{(t\ell)^2}+ O\left( \left(\frac{2}{t\ell}\right)^3 \right)\right) \nonumber \\
&& + t\ell \exp\left(-\ell\left(\frac{1}{t\ell}\right)\right) \exp\left(- \frac{\ell}{2(t\ell)^2}+ O\left( \left(\frac{1}{t\ell}\right)^3 \right)\right) \nonumber\\
&& - (t\ell)^2\exp\left(-2\ell\left(\frac{1}{t\ell}\right)\right) \exp\left(- \frac{\ell}{(t\ell)^2}+ O\left( \left(\frac{1}{t\ell}\right)^3 \right)\right) \nonumber\\
\end{eqnarray}
\begin{eqnarray}
&\sim& (t\ell)(t\ell-1)\exp\left(-\ell\left(\frac{2}{t\ell}\right)\right) \left(1-\frac{2\ell}{(t\ell)^2}\right) \nonumber \\
&& + t\ell \exp\left(-\ell\left(\frac{1}{t\ell}\right)\right) \left(1 - \frac{\ell}{2(t\ell)^2} \right) \nonumber\\
&& - (t\ell)^2\exp\left(-2\ell\left(\frac{1}{t\ell}\right)\right) \left(1- \frac{\ell}{(t\ell)^2}\right) \nonumber\\
&\sim& t \ell e^{-\frac{1}{t}} - t \ell e^{-\frac{2}{t}} - \ell e^{-\frac{2}{t}} \nonumber \\
&\sim&\sigma_n^2.\label{variance}
\end{eqnarray}
Under the assumption that $t=o(\ell^{1/5})$, we can verify that $\frac{\mu_n^2}{\sigma_n^6} \to 0$ as $n \to \infty$.
We note $\mu_n/\sigma_n$ is of the order $\sqrt{t^3\ell}$. Hence for $k$ in the range specified in Lemma~\ref{normal}, we have
$k^2=O(t^3\ell)$. Therefore
\begin{eqnarray*}
(t\ell)_k&=&(t\ell)^k\prod_{j=1}^{k-1}\left(1-\frac{j}{t\ell}\right)\\
&=&(t\ell)^k\exp\left(\sum_{j=1}^{k-1}\ln\left(1-\frac{j}{t\ell}\right)\right)\\
&=&(t\ell)^k\exp\left(-\sum_{j=1}^{k-1}\left(\frac{j}{t\ell}+\frac{j^2}{2t^2\ell^2}+O\left(\frac{j^3}{t^3\ell^3}\right)\right)\right)\\
&\sim&(t\ell)^k\exp\left(-\frac{k^2}{2t\ell}\right),\\
\left(\frac{t\ell-k}{t\ell}\right)^{\ell}&=&\exp\left(\ell\ln\left(1-\frac{k}{t\ell}\right)\right)\\
&=&\exp\left(-\ell\left(\frac{k}{t\ell}+\frac{k^2}{2t^2\ell^2}+O\left(\frac{k^3}{t^3\ell^3}\right)\right)\right)\\
&\sim&\exp\left(-\frac{k}{t}-\frac{k^2}{2t^2\ell}\right).
\end{eqnarray*}
It follows from (\ref{factorial}) and  (\ref{mean}) that
$$\ex((Y_n)_k)\sim (t\ell)^k\exp\left(-\frac{k}{t}-\frac{k^2}{2t\ell}-\frac{k^2}{2t^2\ell}\right)
\sim \mu_n^k\exp\left(-\frac{k^2(t+1)}{2t^2\ell}\right).
$$
Now the result follows from Lemma~\ref{normal} and the estimation
$$s_n=\frac{\sigma_n^2-\mu_n}{\mu_n^2}=\frac{\ex(Y_n(Y_n-1))}{\mu_n^2}-1
=-\frac{1}{t\ell}-\frac{1}{t^2\ell}+O\left(\frac{1}{t^2\ell^2}\right) > - \mu_n^{-1},$$
as $\frac{t}{\ell} \to 0$ when $\ell \to \infty$.
\qed

The following corollary follows immediately from Theorem~\ref{Thm1}, by noting $s=t=1$ and $X_n=n-Y_n$.
\begin{corollary}\label{Small}
Under the assumption of Theorem~\ref{Thm1} and $X_n= n-Y_n$, we have
$$\pr(X_n=h)=\frac{1}{n^n}{n\choose h}\sum_{j=1}^h(-1)^{h-j}{h\choose j}j^n.
$$
In particular,
$$\pr(X_n=n)=\frac{n!}{n^n},
$$
and for  $h=o(n)$,
$$\pr(X_n=h)\sim \frac{1}{h!}n^h\left(\frac{h}{n}\right)^{n}.
$$
\end{corollary}

Consider $\ell=n=q-1$ and $t=1$.  Because an $r$-th order cyclotomic mapping polynomial $f^{r}_{a_0, \ldots, a_{q-2}}(x)$ always maps $0$ to $0$,  Corollary~\ref{Small}
implies

\begin{corollary}\label{Smallwithoutzero} Let $g(x)$ be any random $r$-th order cyclotomic mapping polynomial $f^{r}_{a_0, \ldots, a_{q-2}}(x)$ over finite field $\fq$ such that none of $a_i$'s is zero. Then the probability of $\mid V_f \mid = h+1$ is given by
$$\pr(X_{q-1}=h)= \frac{1}{(q-1)^{q-1}}{q-1\choose h}\sum_{j=1}^h(-1)^{h-j}{h\choose j}j^{q-1}.$$
In particular, the probability of such a random cyclotomic mapping polynomial $f^{r}_{a_0, \ldots, a_{q-2}}(x)$ is a permutation polynomial is
$$\pr(X_{q-1}=q-1)=\frac{(q-1)!}{(q-1)^{q-1}},
$$
and for $h=o(q)$, the probability of such a random cyclotomic mapping polynomial $f^{r}_{a_0, \ldots, a_{q-2}}(x)$ with a value set size $h+1$ is
$$\pr(X_{q-1}=h)\sim \frac{1}{h!}(q-1)^{h}\left(\frac{h}{q-1}\right)^{q-1}.
$$
\end{corollary}

\section{Sizes of value sets of random polynomials }\label{random polynomials}

Recall $q-1 = \ell s$ and $r$ is a positive integer  such that $(r, s) = t$.
In this section we consider the value set size for any random $r$-th order cyclotomic mapping polynomial $f^{r}_{a_0, \ldots, a_{\ell-1}}(x)$ with index $\ell$ over finite field $\fq$. Namely, $a_0, \ldots, a_{\ell-1}$ are independently chosen at random from $\F_q$. 
 We note that the value set problem for any random polynomial with degree at most $q-1$ is in fact the value set problem for a random $r$-th order cyclotomic mapping polynomial with index $\ell=q-1$.

As discussed in Section~\ref{Methodology}, we are interested in the distribution of $X_n=|\cup_{j=1}^{\ell}A_j|$.
However, we need to include the element 0 in our analysis similar to those in Section~\ref{random mappings}. Because each $a_i$ in (\ref{CycloMappingDef}) is chosen independently at random from $\fq$, this leads to the following model.

Let us define $D_0=\{0\}$ and $D_j = T_{j-1}$ for $1\leq j \leq t\ell$.
The distribution of each random set $A_i$ is $\pr(A_i=D_0)=\frac{1}{q}$,  and
$$\pr(A_i=D_j)=\frac{1}{t\ell}\left(1-\frac{1}{q}\right)=\frac{s}{tq},~j=1,\ldots,\ell.$$
As in the case discussed earlier, we define the events $B_j= \cap_{i=1}^\ell \{A_i\ne D_j\}$, $0\le j\le t\ell$.
 Let $Y_{t\ell}$ be the number of $B_1,B_2,\ldots, B_{t\ell}$ (note $B_0$ is excluded) which occur.
Then $\pr(Y_{t\ell}=k)$ is the probability that exactly $k$ of $t\ell$ cyclotomic sets $T_j$'s of index $t \ell$ are not in the value set of $g$. Then we have
\begin{lemma}\label{randomExpected}
 Let $q-1 = \ell s$ and $r$ be a positive integer  such that $(r, s) = t$.  Let $g(x)$ be any random $r$-th order cyclotomic mapping polynomial $f^{r}_{a_0, \ldots, a_{\ell-1}}(x)$ with index $\ell$ over finite field $\fq$. Let $Y_{t\ell}$ be the number of cyclotomic sets of index $t \ell$ not contained in the value set of $g$. Then
\begin{eqnarray*}
\ex((Y_{t\ell})_k)&=&(t\ell)_k\left(1-\frac{sk}{tq} \right)^{\ell},\\
\pr(Y_{t\ell}=k)&=&{t\ell\choose k}  \sum_{j=0}^{t\ell - k} (-1)^{t\ell -j-k} {t\ell-k\choose j}\left(\frac{1}{q} + \frac{sj}{tq}\right)^{\ell}.
\end{eqnarray*}
\end{lemma}
\proof Let $\CN = \{1, 2, \ldots, t\ell\}$. Define
$$
S_k=\sum_{J\subset \CN,|J|=k}\pr(\cap_{j\in J}B_j).
$$
Because $A_1, \ldots, A_{\ell}$ are mutually independent, we have
$$
S_k={t\ell\choose k}\pr(\cap_{j=1}^k B_j)={t\ell\choose k}\left(\frac{1}{q}+\frac{s(t\ell-k)}{tq}\right)^{\ell}={t\ell\choose k}\left(1-\frac{sk}{tq} \right)^{\ell}.
$$

Using Equation~(\ref{pk}), we obtain

\begin{eqnarray*}
\pr(Y_{t\ell}=k)&=&\sum_{h=k}^{t\ell} (-1)^{h-k} { h \choose k} S_h \\
&=& \sum_{h=k}^{t\ell} (-1)^{h-k} { h \choose k}  {t\ell\choose h}\left(\frac{1}{q}+\frac{s(t\ell-h)}{tq}\right)^{\ell} \\
&=&  {t\ell\choose k} \sum_{h=k}^{t\ell} (-1)^{h-k} {t\ell-k\choose h-k} \left(\frac{1}{q}+\frac{s(t\ell-h)}{tq}\right)^{\ell} \\
&=& {t\ell\choose k}  \sum_{j=0}^{t\ell - k} (-1)^{t\ell -j-k} {t\ell-k\choose j}\left(\frac{1}{q} + \frac{sj}{tq}\right)^{\ell}.
\end{eqnarray*}

 \qed

Using the above lemma, we can obtain the distribution of $X_{t\ell} = |V_g|$.
\begin{theorem}\label{valuesetexclude0}
 Let $q-1 = \ell s$ and $r$ be a positive integer  such that $(r, s) = t$.  Let $f(x)$ be any random $r$-th order cyclotomic mapping polynomial $f^{r}_{a_0, \ldots, a_{\ell-1}}(x)$ with index $\ell$ over $\fq$.
 Then
 $$\pr(X_{t\ell}=1+ks/t) = {t\ell\choose k}  \sum_{j=0}^{k} (-1)^{k -j} {k\choose j}\left(\frac{1}{q} + \frac{sj}{tq}\right)^{\ell}.
  $$
\end{theorem}
\proof  Indeed,
\begin{eqnarray*}
\pr(X_{t\ell}=1+ks/t)&=&\pr(\{Y_{t\ell}=t\ell-k\})\\
&=& {t\ell\choose k}  \sum_{j=0}^{k} (-1)^{k -j} {k\choose j}\left(\frac{1}{q} + \frac{sj}{tq}\right)^{\ell}.
\end{eqnarray*}

We now obtain an application of the above results on random polynomials with degree at most $q-1$. It is sufficient to study the distribution of the subclass of random polynomials $g(x)$ such that $g(0) =0$ with degree $\leq q-1$, because any polynomial $f(x)$ such that $f(0) =b$ can be written as $g(x) +b$ with $g(0) =0$ and vice versa.  So we can view any random polynomial $g$ of degree at most $q-1$ with $g(0)=0$ and index $\ell$ as an $r$-th order cyclotomic mapping polynomial with the least index $\ell$, which is therefore a cyclotomic mapping polynomial with index $q-1$. Because there are $q^{q-1}$ such polynomials with $g(0) =0$, they correspond to all the cyclotomic mapping polynomials with index $q-1$. So a random polynomial with degree less than or equal to $q-1$ and $g(0)=0$ is a random $r$-th order cyclotomic mapping polynomial with index $\ell=q-1$ for any $r \geq 1$.  Therefore,  Lemma~\ref{randomExpected} implies that any random polynomial with degree $q-1$ has expected value set size $(1-\frac{1}{q})^{q-1}\sim \frac{q}{e}$. This verifies William's result \cite{Williams} saying that almost all the polynomials of degree $q-1$ is a general polynomial.  Moreover, applying Theorem~\ref{valuesetexclude0} to the case $\ell =q-1$ (hence $s=t=1$),  we obtain exact probability distribution of the value set size for a random polynomial over finite field $\fq$ in Corollary~\ref{Smallwithzero}.
Moreover, we can drive Corollary~\ref{randomnormaldistribution} that the distribution of $(Y_{t\ell}-u_{t\ell})/\sigma_{t\ell}$ tends to the standard normal distribution as well, as $\ell \to \infty$ and $t = o(\ell^{1/5})$,  following the same arguments as in the proof of Theorem~\ref{normaldistribution}.


\section{Size of the Union of Random Sets}\label{union of random sets}

Recall $\CN=\{1, \ldots, n\}$.
In Section~\ref{random mappings}, we considered
$|\cup_{i=1}^{\ell} A_i|$  where each $A_i$ a random cyclotomic cosets of index $\ell$. In Section~\ref{random polynomials},
we allowed the possibility that $A_i$ could be $\{0\}$, a subset with a different size from cyclotomic cosets. Earlier,   Barot and Pe\~{n}a \cite{BP} considered the probability distribution of
$|\cup_{i=1}^{\ell} A_i|$ where each $A_i$ is chosen independently and uniformly at random from $\mathcal{P}_{m}$ where $\mathcal{P}_m$ is the set of all $m$-subsets of $\CN$.
In this section, we consider $|\cup_{i=1}^{\ell} A_i|$ such that $m_i$'s can be distinct. Let $X_n=|\cup_{i=1}^{\ell} A_i|$ and $Y_n=n-X_n$. In the following we establish the distributions of  $X_n$ and $Y_n$ and we show that the $Y_n$  is asymptotically normal. This generalizes the following result by Barot and Pe\~{n}a \cite{BP} because we allow $m_i$'s to be distinct.

\begin{theorem}[Barot and Pe\~{n}a, 2001]
Let $\CN=\{1, \ldots, n\}$.  Let $X_n=|\cup_{i=1}^{\ell} A_i|$ where each $A_i$ is chosen independently and uniformly at random from $\mathcal{P}_{m}$ where $\mathcal{P}_m$ is the set of all $m$-subsets of $\CN$ and $Y_n=n-X_n$.
We have
\begin{eqnarray*}
\pr(X_n=i) &=&\frac{{n\choose i}}{{n\choose m}^\ell}\sum_{j=0}^{i-m}(-1)^{j}{i\choose j}{i-j\choose m}^\ell,\\
\ex(X_n) &=&n\left(1-\left(1-\frac{m}{n}\right)^\ell\right),\\
\var(X_n) &=&n(n-1)\left(1-\frac{m}{n}\right)^\ell \left(1-\frac{m}{n-1}\right)^\ell -\left(\ex(X_n)\right)^2+\ex(X_n).
\end{eqnarray*}
\end{theorem}

First of all, we extend the above results to general $m_i$'s.

\begin{lemma} Let $\CN=\{1, \ldots, n\}$ and  $\mathcal{P}_{m_i}$ be the set of all $m_i$-subsets of $\CN$
for $1\leq i\leq \ell$. Let $X_n$ be the random variable for the size of $\cup_{i=1}^{\ell} A_i$ where $A_i$ is a random set chosen independently and uniformly from $\mathcal{P}_{m_i}$ and $Y_n=n-X_n$. We have
\begin{eqnarray}
\ex((Y_n)_k)&=&(n)_k\prod_{j=1}^\ell \frac{(n-m_j)_k}{(n)_k},\\
\pr(Y_n=k)&=&\sum_{h=k}^n(-1)^{h-k}{h \choose k}{n\choose h}\prod_{j=1}^\ell \frac{(n- m_j)_h}{(n)_h},\\
\pr(X_n=i)&=&{n\choose i}\sum_{h=0}^{i}(-1)^h{i\choose h}\prod_{j=1}^\ell \frac{ {i-h \choose m_j}}{{n \choose m_j}}.
\end{eqnarray}
\end{lemma}
\proof  For each $j\in \CN$, let $B_j$ denote the event that $j\notin \cap_{i=1}^\ell A_i$. We note that $\pr(Y_n=k)$ is the probability that exactly $k$ of the $B_j$'s occur.
Since $A_1,A_2,\ldots,A_\ell$ are mutually independent, we have
\begin{eqnarray*}
S_k&=&{n\choose k}\pr(B_1\cap B_2\cap \cdots \cap B_k)\\
&=&{n\choose k}\prod_{j=1}^\ell \pr(1\notin A_j,2\notin A_j,\ldots,k\notin A_j)\\
&=&{n\choose k}\prod_{j=1}^\ell \frac{(n-m_j)_k}{(n)_k},
\end{eqnarray*}
and hence
$$\ex((Y_n)_k)=k!S_k=(n)_k\prod_{j=1}^\ell \frac{(n-m_j)_k}{(n)_k}.
$$
Now from (\ref{pk}) we obtain
$\pr(Y_n=k) =\sum_{h=k}^n(-1)^{h-k}{h \choose k}{n\choose h}\prod_{j=1}^\ell \frac{(n- m_j)_h}{(n)_h}$.
Finally,
\begin{eqnarray*}
\pr(X_n=i)&=&\pr(Y_n = n-i) \\
&=& \sum_{h=n-i}^n(-1)^{h-n+i}{h \choose n-i}{n\choose h}\prod_{j=1}^\ell \frac{(n- m_j)_h}{(n)_h}\\
&=&{n\choose i}\sum_{h=0}^{i}(-1)^h{i\choose h}\prod_{j=1}^\ell \frac{ {i-h \choose m_j}}{{n \choose m_j}},
\end{eqnarray*}
where the last equality holds after the substitution $h: = h-n+i$. \qed.

Next we show that $Y_n$ is asymptotically normal. The condition on $u_i$ in the following theorem can be relaxed considerably, but we
use the current form for the sake of simplicity.
\begin{theorem}
Let $u_j=m_j/n$. Suppose $\ell \to \infty$, $u_j=O(1/\ell)$ uniformly for all $1\le j\le \ell$, and $\sum_{i=1}^\ell u_i>c$ for some positive constant $c$.
Then $Y_n$ is asymptotically normal with mean and variance, respectively, equal to
\begin{eqnarray*}
\mu_n&=&n\prod_{i=1}^\ell (1-u_i),\\
\sigma_n^2&=&n\left(1-\left(1+\sum_{i=1}^\ell u_i\right)\prod_{i=1}^\ell (1-u_i)\right)\prod_{i=1}^\ell (1-u_i).
\end{eqnarray*}
\end{theorem}
\proof For $k=O(\sqrt{n})$,
We have
\begin{eqnarray*}
\ex((Y_n)_k)&=&(n)_k\prod_{i=1}^\ell \frac{(n-m_i)_k}{(n)_k}\\
&=&n^k\exp\left(-\frac{k(k-1)}{2n}+O\left(\frac{k^3}{n^2}\right)\right)\prod_{i=1}^\ell \prod_{j=0}^{k-1}\left(1-u_i-\frac{ju_i}{n-j}\right)\\
&=& \mu_n^k\exp\left(-\frac{k(k-1)}{2n}+O\left(\frac{k^3}{n^2}\right)\right)\prod_{i=1}^\ell \exp\left(\sum_{j=0}^{k-1}
\ln\left(1-\frac{ju_i}{(1-u_i)(n-j)}\right)\right)\\
&=&\mu_n^k\exp\left(-\frac{k(k-1)}{2n}+O\left(\frac{k^3}{n^2}\right)-\sum_{i=1}^\ell \sum_{j=0}^{k-1}\frac{ju_i}{(1-u_i)(n-j)}\right) \\
&=& \mu_n^k\exp\left(-\frac{k(k-1)}{2n}-\sum_{i=1}^\ell \sum_{j=0}^{k-1}\frac{ju_i}{(1-u_i)n}+O\left(\frac{k^3}{n^2}\right)\right) ~~~~ (because ~j \leq k ~\sim \sqrt{n})\\
&=& \mu_n^k\exp\left(-\frac{k(k-1)}{2n}\left(1+\sum_{i=1}^\ell u_i\right)
+O\left(\frac{k^3}{n^2}+\frac{k^2}{n\ell}\right)\right). ~~~(note ~ u_j = O(1/\ell))
\end{eqnarray*}
In particular we have
$$
\ex((Y_n)_2)=\mu_n^2\exp\left(-\frac{1}{n}\left(1+\sum_{i=1}^\ell u_i\right)+O\left(\frac{1}{n^2}+\frac{1}{n\ell}\right)\right),
$$
and for $k$ of the order $\sqrt{n}$, 
$$
\ex((Y_n)_k)\sim \mu_n^k\exp\left(-\frac{k^2}{2n}\left(1+\sum_{i=1}^\ell u_i\right)\right).
$$
It follows that
\begin{eqnarray*}
s_n&=&\frac{\sigma_n^2-\mu_n}{\mu_n^2}\\
&=&\frac{\ex(Y_n(Y_n-1))}{\mu_n^2}-1\\
&=&\exp\left(-\frac{1}{n}\left(1+\sum_{i=1}^\ell u_i\right)+O\left(\frac{1}{n^2}+\frac{1}{n\ell}\right)\right)-1\\
&=&-\frac{1}{n}\left(1+\sum_{i=1}^\ell u_i\right)+O\left(\frac{1}{n^2}+\frac{1}{n\ell}\right),
\end{eqnarray*}
and
$$
\ex((Y_n)_k)\sim \mu_n^k\exp\left(\frac{k^2s_n}{2}\right).
$$
Now the theorem follows from Lemma~\ref{normal}. \qed


\vskip 20pt


\vskip 20pt

\baselineskip 12pt \frenchspacing
\begin{thebibliography}{99}


\bibitem{AGW:09} A. Akbary, D. Ghioca, and Q. Wang, On permutation polynomials of prescribed shape,
{\it Finite Fields Appl.} 15 (2009),  195-206.




\bibitem{BW} B. J. Birch and H. P. F. Swinnerton-Dyer, Note on a problem of Chowla, Acta Arith. 5 (1959), 417-423.

\bibitem{BC} H. Borges and R. Conceicao, On the characterization of minimal value set polynomials, {\it J. Number Theory} 133 (2013), 2021-2035.
\bibitem{BP} M. Barot and J. Pe\~{n}a, Estimating the size of a union of random subsets of fixed cardinality,
{\em Elemente der Mathematik} 56 (2001), no. 4, 163-169.
\bibitem{Bo}  B. Bollob\'{a}s, Random Graphs, Academic Press, 1985.

\bibitem{Carlitzetal:61}
L. Carlitz,  D. J. Lewis, W. H. Mills, and E. G.  Straus, {Polynomials over finite fields with minimal value sets},
 {\it Mathematika} 8 (1961), 121-130.




 \bibitem{CHW} Q. Cheng, J. Hill and D. Wan,  Counting value sets: algorithms and complexity, Tenth Algorithmic Number Theory Symposium ANTS-X, 2012, University of California at San Deigo.



\bibitem{Cohen:70} S. D. Cohen, The distribution of polynomials over finite fields,  {\it Acta Arith.} 17 (1970), 255-271.


\bibitem{Da} F.N. David, {\em Biometrika} 37 (1950), 97-110.

\bibitem{DasMul} P. Das and G. L. Mullen, Value sets of polynomials over finite fields,  in {\it Finite Fields with Applications in Coding Theory, Cryptography and Related Areas},  G.L. Mullen, H. Stichtenoth, and H. Tapia-Recillas, Eds., Springer, 2002, 80-85.




\bibitem{Evans:92} A. B. Evans, Orthomorphism Graphs of Groups, Lecture Notes in Mathematics, Vol. 1535, Springer, Berlin, 1992.



\bibitem{GW} Z. Gao and N.C. Wormald, Asymptotic normality determined by high
moments, and submap counts of random maps, {\em Probab. Theory Relat. Fields} 130 (2004), 368-376.



\bibitem{GomezMadden:88} J. Gomez-Calderon and D. J. Madden,
Polynomials with small value set over finite fields,  {\it J. Number Theory}  28 (1988), no. 2, 167-188.

\bibitem{GurWan} R. Guralnick and D. Wan, Bounds for fixed point
free elements in a transitive group and applications to curves over finite fields,
{\it Israel J. Math.} 101 (1997),  255-287.

\bibitem{Mills:64} W. H. Mills,  Polynomials with minimal value sets,
{\it Pacific J. Math} 14 (1964), 225-241.










\bibitem{MWW} G. L. Mullen, D. Wan, and Q. Wang, Value sets of polynomial maps over finite fields, {\it  Quart. J. Math.}  64 (2013), no. 4, 1191-1196.

\bibitem{MWW2} G. L. Mullen, D. Wan, and Q. Wang, An index bound on value sets of polynomial maps over finite fields, {\it  Proceedings of Workshop on the Occasion of Harald Niederreiter's 70th Birthday: Applications of Algebra and Number Theory}, June 23-27, 2014.


\bibitem{NW:05} H. Niederreiter and A. Winterhof,
Cyclotomic $\mathcal{R}$-orthomorphisms of finite fields,
{\it Discrete Math.}  295 (2005), 161-171.
 
\bibitem{Tu} G. Turnwald, A new criterion for permutation polynomials,  {\it Finite Fields  Appl.} 1 (1995), 64-82.


\bibitem{W1}D. Wan, A $p$-adic lifting lemma and its applications to permutation
polynomials, Lecture Notes in Pure and Appl. Math.,  Marcel Dekker,
New York, Vol. 141, 1992, 209-216.


\bibitem{WSC}D. Wan, P. J. S. Shiue and C. S. Chen, Value sets of polynomials over
finite fields, {\it Proc. Amer. Math. Soc}. 119 (1993), 711-717.




\bibitem{Wang}
Q. Wang, \emph{Cyclotomic mapping permutation polynomials over finite fields}, Sequences, Subsequences, and
Consequences (International Workshop, SSC 2007, Los Angeles, CA, USA, May 31 - June 2, 2007), 119-128,
 Lecture Notes in Comput. Sci.  Vol. 4893, Springer, Berlin, 2007.


\bibitem{Wang2}
Q. Wang,  Cyclotomy and permutation polynomials of large indices, {\it Finite Fields Appl.}  22 (2013), 57-69.




\bibitem{Williams}
K. S. Williams, On general polynomials, {\it Canad. Math. Bull.} 10 (1967), no. 4, 579-583.








\end{thebibliography}
\end{document}